\documentclass[12pt,a4paper]{article}
\usepackage{amsmath}
\usepackage{amsfonts,amssymb,amscd,amsthm}
\usepackage{url}

\theoremstyle{plain}
\newtheorem{theorem}{Theorem}

\newtheorem{lemma}[theorem]{Lemma}

\theoremstyle{definition}
\newtheorem{definition}[theorem]{Definition}
\theoremstyle{remark}
\newtheorem{remark}[theorem]{Remark}

\def\CC{\mathbb{C}}

\def\cB{\mathcal{B}}

\def\bcdot{\,\boldsymbol\cdot\,}
\def\lra{\longrightarrow}

\def\si{\sigma}
\def\ph{\varphi}
\def\Ga{\Gamma}
\def\La{\Lambda}
\def\om{\omega}

\DeclareMathOperator\dist{dist}

\DeclareMathOperator\tr{tr} 
\DeclareMathOperator\ind{ind}
\DeclareMathOperator\Vect{Vect}
\DeclareMathOperator\End{End}
\DeclareMathOperator\id{id}

\newcommand{\op}[1]{\operatorname{#1}}

\def\ov{\overline}
\def\ovs{\overset}

\newcommand{\BL}{\biggl}
\newcommand{\BR}{\biggr}
\newcommand{\bl}{\bigl}
\newcommand{\br}{\bigr}

\newcommand{\abs}[1]{\lvert#1\rvert}

\def\Td{\operatorname{Td}}

\def\ch{\operatorname{ch}}

\title{On Elliptic Differential Operators with Shifts \\[2pt]
       II. The Cohomological Index Formula}
\author{V.~E.~Nazaikinskii, A.~Yu.~Savin, and B.~Yu.~Sternin}
\date{}

\begin{document}
\maketitle

\subsection*{Introduction}

This paper is a continuation of \cite{0706.3511}, where we
have studied a general class of (pseudo)\-diffe\-ren\-tial
operators with nonlocal coefficients, referred to as
operators with shifts, and obtained a local index formula
(i.e., a formula expressing the index as the integral of a
differential form explicitly determined by the principal
symbol of the operator) for matrix elliptic operators of
this kind. In the present paper we finish the business by
establishing a cohomological index formula of
Atiyah--Singer type for elliptic differential operators
with shifts acting between section spaces of arbitrary
vector bundles. The key step is the construction of closed
graded traces on certain differential algebras over the
symbol algebra for this class of operators.

We do not formally assume the reader to be familiar with
\cite{0706.3511} as far as definitions are concerned but
freely use the results obtained there. We also do not
reproduce the discussion of general motivations for this
research, which, as well as the bibliography, can be found
in \cite{0706.3511}.

\paragraph{Acknowledgements.}
The research was supported in part by RFBR grants
nos.~05-01-00982 and~06-01-00098 and DFG grant 436 RUS
113/849/0-1\circledR ``$K$-theory and noncommutative
geometry of stratified manifolds."

The authors thank Professor Schrohe and Leibniz
Universit\"at Hannover for kind hospitality.

\section{Elliptic operators with shifts}

\subsection{Pseudodifferential operators with shifts}

\paragraph{The group $\Ga$.}

Let $M$ be a compact oriented Riemannian manifold without
boundary, and let $\Ga$ be a countable dense subgroup of a
Lie group $\ov\Ga$ of orientation-preserving isometries of
$M$. The natural action of $\ov\Ga$ on functions on $M$
will be denoted by $T$, so that
\begin{equation*}
    [T_gu](x)=u(g^{-1}(x)),\qquad x\in M.
\end{equation*}
We assume that $\Ga$ satisfies the following two
conditions:
\begin{enumerate}
    \item (\textit{Polynomial growth}.) The group $\Ga$
    is finitely generated, and the number of distinct
    elements of $\Ga$ representable by words of length $\le k$ in
    some finite system of generators grows at most polynomially in $k$.

    In what follows, we fix some system of generators and denote by
    $\abs{g}$ the minimum length of words representing $g\in\Ga$.
    \item (\textit{Diophantine property}.) Let $\op{fix}(g)$
be the set of fixed points of $g\in\Ga$. The estimate
\begin{equation*}
    \dist(g(x),x)\ge C\abs{g}^{-N}\dist(x,\op{fix}(g))
\end{equation*}
holds for some $N,C>0$ and for all $x\in M$ and $g\in\Ga$.
Here $\dist(x,\op{fix}(g))$ is the Riemannian distance
between $x$ and the set $\op{fix}(g)$, and by convention we
set $\dist(x,\op{fix}(g))=1$ if $\op{fix}(g)$ is empty.
\end{enumerate}

\paragraph{Matrix operators.}
Matrix pseudodifferential operators with shifts, \i.e.,
$\Psi$DO with shifts acting on vector functions on $M$, can
be described as follows. (For more detail, see
\cite{0706.3511}, where also further bibliographical
references can be found.) A matrix $\Psi$DO of order $m$
with shifts has the form
\begin{equation}\label{shifts1}
D=\sum_{g\in\Gamma} T_g D_g,
\end{equation}
where $D_g$ is a classical $\Psi$DO of order $m$ on $M$ and
the operators $D_g$ rapidly decay as $\abs{g}\to\infty$ in
the natural Fr\'echet topology on the set of $m$th-order
$\Psi$DO.

\paragraph{Operators on sections of vector bundles.}
Pseudodifferential operators with shifts acting on sections
of vector bundles are an easy generalization of matrix
operators. To define them, one should localize into
neighborhoods where the bundles are trivial. The only
difference with the case of pseudodifferential operators
without shifts is that our operators are no longer local,
so we cannot localize into a neighborhood of the diagonal;
hence two neighborhoods, instead of one, in the subsequent
argument. Let $E$ and $F$ be finite-dimensional complex
vector bundles on $M$. A linear operator
\begin{equation}\label{shi2}
 D\colon C^\infty(M,E)\longrightarrow C^\infty(M,F)
\end{equation}
is called an $m$th-order $\Psi$DO with shifts if for any
trivializations of $E$ and $F$ over some neighborhoods
$U_E,U_F\subset M$, respectively, and any functions
$\varphi\in C^\infty_0(U_E)$ and $\psi\in C_0^\infty(U_F)$
the operator $\psi D\ph$ is an $m$th-order matrix $\Psi$DO
with shifts of the form~\eqref{shifts1}.

We point out that no action of $\Ga$ on the bundles $E$ and
$F$ is needed in this definition.

The linear space of $m$th-order pseudodifferential
operators \eqref{shi2} with shifts will be denoted by
$\Psi^m(E,F)_\Gamma$. If $E,F$ and $H$ are three vector
bundles on $M$, then the multiplication of operators
induces a well-defined bilinear mapping
\begin{equation*}
 \Psi^m(E,F)_\Ga\times \Psi^{m'}(F,H)_\Ga\lra
\Psi^{m+m'}(E,H)_\Ga.
\end{equation*}
Just as for matrix operators, one readily proves that an
$m$th-order $\Psi$DO with shifts is a continuous operator
of order $m$ in the Sobolev spaces of sections of $E$ and
$F$.

\subsection{Symbol, ellipticity, and Fredholm property}

\paragraph{Symbol: the matrix case.}

First, let us recall what happens in case the bundles $E$
and $F$ are trivial.

For the $n\times n'$ matrix operator \eqref{shifts1}, the
symbol is defined by the formula
\begin{equation}\label{symbol2}
\sigma(D)=\sum_{g\in \Ga}T_{\partial g}\sigma(D_g)\colon
L^2(S^*M,\CC^n)\lra L^2(S^*M,\CC^{n'}),
\end{equation}
where the codifferential
\begin{equation*}
 \partial g\colon S^*M\to S^*M
\end{equation*}
is the map induced by $g$ (it acts as $g$ along the base
and as $((dg)^*)^{-1}$ in the fibers of $S^*M$).

\paragraph{Symbol: the general case.}

If the operator \eqref{shi2} is a usual pseudodifferential
operator, then its symbol is a bundle homomorphism
$\pi^*E\to\pi^*F$, where $\pi\colon S^*M\to M$ is the
natural projection. For pseudodifferential operators with
shifts, which are highly nonlocal, this is no longer the
case, and their symbols are defined as homomorphisms of
section spaces of the bundles $\pi^*E$ and $\pi^*F$ rather
than of the bundles themselves.
\begin{definition}
The \textit{symbol} of the operator \eqref{shi2} is the operator
\begin{equation}\label{symbol1}
\sigma(D)\colon L^2(S^*M,\pi^*E)\lra L^2(S^*M,\pi^*F)
\end{equation}
such that for any trivializations of $E$ and $F$ over some
neighborhoods $U_E,U_F\subset M$, respectively, and any
functions $\varphi\in C^\infty_0(U_E)$ and $\psi\in
C_0^\infty(U_F)$ the operator $\psi\sigma(D)\varphi$ is the
symbol of the operator  $\psi D\ph$.
\end{definition}

One can readily verify that the symbol of a $\Psi$DO with
shifts is well defined. The space of symbols of $\Psi$DO
with shifts acting between section spaces of vector bundles
$E$ and $F$ will be denoted by
$C^\infty(S^*M,\op{Hom}(E,F))_\Gamma$. For $E=F$, we use
the notation $C^\infty(S^*M,\End(E))_\Gamma$, and for
scalar symbols write $C^\infty(S^*M)_\Gamma$, just as in
the first part of the paper. A generalization of the
argument given there shows that
$C^\infty(S^*M,\End(E))_\Gamma$ is a local subalgebra of
the $C^*$-algebra $\cB L^2(S^*M,E)$. Hence if a symbol
\begin{equation*}
 \si \in C^\infty(S^*M,\op{Hom}(E,F))_\Gamma
\end{equation*}
is invertible (as an operator in $L^2$), then one
necessarily has
\begin{equation*}
 \si^{-1} \in C^\infty(S^*M,\op{Hom}(F,E))_\Gamma.
\end{equation*}

\begin{definition}
An operator $D\in\Psi^m(E,F)_\Gamma$ is said to be
\textit{elliptic} if its symbol $\si(D)$ is invertible.
\end{definition}

As usual, one has the finiteness theorem.

\begin{theorem}[the finiteness theorem]
An operator $D\in\Psi^m(E,F)_\Gamma$ is Fredholm if and
only if its symbol is invertible.
\end{theorem}

\section{The index theorem}

In this section we obtain a cohomological index formula for
elliptic operators $D\in\Psi^m(E,F)_\Gamma$. First, we
shall introduce the elements that occur in this formula.

\subsection{Some objects associated with the group $\Ga$}

We represent the group $\Ga$ as the disjoint union
\begin{equation*}
    \Ga=\bigsqcup_{g_0}\langle g_0\rangle
\end{equation*}
of conjugacy classes and arbitrarily fix an element, $g_0$,
in each conjugacy class $\langle g_0\rangle$. In what
follows, the symbol $g_0$ is invariably used to denote this
fixed representative. By $C_{g_0}$ we denote the
centralizer of $g_0$ in $\ov\Ga$:
\begin{equation*}
    C_{g_0}=\{h\in\ov\Ga\colon hg_0h^{-1}=g_0\}.
\end{equation*}
This is a closed Lie subgroup of $\ov\Ga$. For each
$g\in\langle g_0\rangle$, consider the set $\ov\Ga_{g_0,g}$
of elements $h\in\ov\Ga$ conjugating $g_0$ with $g$, that
is, satisfying
\begin{equation*}
 hg_0h^{-1}=g.
\end{equation*}
Clearly, $\ov\Ga_{g_0,g}$ is a left coset of $C_{g_0}$ in
$\ov\Ga$ and, as such, has a well-defined normalized Haar
measure $dh$ induced by that on $C_{g_0}$.

If the group $\Ga$ acts on a compact manifold $X$, then by
$X_g$ we denote the set of fixed points of an element
$g\in\Ga$. This is a $C^\infty$ submanifold of $X$
consisting of finitely many components (possibly of various
dimensions).

\subsection{The Todd class}

The Todd class $\Td(TM\otimes \mathbb{C};\Gamma)$ of the
complexified tangent bundle of $M$ with respect to the
action of $\Ga$ is an element of the group
$\prod_{g_0}H^{ev}(M_{g_0},\mathbb{C})$. (The product is
taken over representatives of all conjugacy classes in
$\Ga$.) The $g_0$th component of the Todd class is defined
by the formula
\begin{equation}\label{Todd}
 \Td(TM\otimes \mathbb{C};\Gamma)(g_0) =\frac
{\Td(T^*M_{g_0}\otimes
\mathbb{C})}{\ch\lambda_{-1}(NM_{g_0}\otimes
\mathbb{C})(g_0)}\in H^{ev}(M_{g_0},\mathbb{C})
\end{equation}
(This form was apparently first introduced by Atiyah and
Singer in~\cite{AtSi3}; following Baum and Connes
\cite{BaCo2}, we refer to it as the ``Todd class.'')

Let us make some explanations concerning this formula. The
numerator is the usual Todd class of the complexified
tangent bundle of $M_{g_0}$. Next, $\lambda_{-1}(NM_{g_0})$
is the (virtual) vector bundle
\begin{equation*}
 \lambda_{-1}(NM_{g_0})
 = \La^{even}(NM_{g_0})-\La^{odd}(NM_{g_0})
\end{equation*}
composed of the exterior powers of $NM_{g_0}$, and
$\ch\lambda_{-1}(NM_{g_0}\otimes \mathbb{C})(g_0)$ is the
Chern character of the bundle
$\lambda_{-1}(NM_{g_0})\otimes \mathbb{C}$ localized at the
element $g_0$. Recall that it is defined as follows. Since
the mapping $g_0$ preserves the metric, it follows that the
restriction of the differential $dg_0$ to the normal bundle
$NM_{g_0}$ is a well-defined automorphism of this bundle.
Let $\Omega$ be the curvature form of some $dg_0$-invariant
connection on $\lambda_{-1}(NM_{g_0})$ (e.g., of the
connection induced by the restriction of the Riemannian
connection on $TM$ to $NM_{g_0}$). The localized Chern
character
\begin{equation*}
 \ch\lambda_{-1}(NM_{g_0}\otimes
\mathbb{C})(g_0)\in H^{ev}(M_{g_0},\mathbb{C})
\end{equation*}
is defined as the cohomology class of the
form
\begin{equation*}
\ch\lambda_{-1}(NM_{g_0}\otimes \mathbb{C})(g_0)
 =\tr\BL( dg_0^* \exp\BL(-\frac 1{2\pi
i}\Omega\BR)\BR).
\end{equation*}
(Here $\tr$ stands for the trace in the fibers of a vector
bundle.)

\subsection{The Chern character of the symbol}

Let the group $\Ga$ act on a compact manifold $X$.

\paragraph{Differential forms and graded traces over the algebra
$C^\infty(X)_\Gamma$.}

Let $E\in\Vect(X)$ be a vector bundle. By
\begin{equation*}
 \Lambda^*(X,\End
E)_\Gamma\subset \mathcal{B}L^2(X,\Lambda^*(X)\otimes E)
\end{equation*}
we denote the subalgebra of elements $A$ of the form
\begin{equation*}
 A=\sum_{g\in\Gamma}\om_g,
\end{equation*}
where the $\om_g$ have the following property: for any two
functions $\psi$ and $\varphi$ with supports in
neighborhoods where $E$ is trivialized, one has
$$
\om_g= T_g a_g,
$$
where $T_g\omega:=(g^*)^{-1}\omega$ and $a_g$ are some
differential forms on $X$ rapidly decaying in the
$C^\infty$ Fr\'echet topology as $\abs{g}\to\infty$.

We define a mapping
\begin{equation}\label{trace1}
\tau\colon\Lambda^*(X,\End E)_\Gamma\lra \bigoplus_{g_0}
\Lambda^*(X_{g_0})
\end{equation}
(the sum is taken over representatives of all conjugacy
classes in $\Ga$) by setting
$$
\tau\BL(\sum_{g\in \Gamma} \omega_g ,g_0\BR) =
\sum_{g\in\langle g_0\rangle}
\int_{\overline{\Gamma}_{g_0,g}}h^*(\omega_g{\bigm|}_{X_g})\,dh.
$$
This is well defined. Indeed, $g|_{X_g}=\op{id}$, and so
the operator $\omega_g$ can be restricted to $X_g$, the
restriction being an $\End E$-valued differential form on
$X_g$. The trace $\tr$ in the last formula is the fiberwise
trace in $\End E$.
\begin{lemma}
The mapping $\tau$ is a graded trace on the algebra
$\Lambda^*(X,\End E)_\Gamma$ in the sense that
\begin{equation*}
 \tau\bl([\omega_1,\omega_2]\br)=0,
\end{equation*}
for all $\omega_1,\omega_2\in \Lambda^*(X,\End E)_\Gamma$,
where $[\bcdot,\bcdot]$ is the supercommutator
$$
[\omega_1,\omega_2]=\omega_1\omega_2-(-1)^{\deg \omega_1 \deg
\omega_2}\omega_2\omega_1.
$$
\end{lemma}

\paragraph{Chern character of projections.}
Now we shall define the \textit{Chern characte}r
$$
\ch\colon K_0(C^\infty(X,\End E)_\Gamma) \longrightarrow
\bigoplus_{g_0} H^{ev}(X_{g_0}, \mathbb{C})
$$
(where $K_0(A)$ is the $K$-group of an operator algebra
$A$). Let $p$ be a projection over the algebra
$C^\infty(X,\End E)_\Gamma$. (To make the subsequent
formulas shorter, we pretend that $p$ is a projection in
the algebra $C^\infty(X,\End E)_\Gamma$ itself rather than
in a matrix algebra over it.) We take some connection
\begin{equation*}
 \nabla_E:\Lambda^*(X,E)\longrightarrow \Lambda^*(X,E)
\end{equation*}
in the bundle $E$ and define a first-order differential
operator with shifts,
\begin{equation}\label{nabbla}
\nabla\colon\Lambda^*(X,E)\longrightarrow \Lambda^*(X,E),
\end{equation}
by the formula
\begin{equation}\label{nabbla1}
    \nabla= p\nabla_E p.
\end{equation}
A straightforward computation shows that the following
assertion is true.
\begin{lemma}\label{lemB}
The operator
$$
\Omega\equiv\nabla^2\colon \Lambda^*(X, E)\lra
\Lambda^*(X,E)
$$
belongs to $\Lambda^2(X,\End E)_\Gamma$.
\end{lemma}
The noncommutative $2$-form $\Omega$ is called the
\emph{curvature form} corresponding to the projection $p$
and the connection $\nabla_E$.

\begin{definition}\label{chern1}
The \emph{Chern character of the class} $[p]\in
K_0(C^\infty(X,\End E)_\Gamma)$ is the cohomology class
\begin{equation*}
    \ch_\Ga [p] \in \bigoplus_{g_0} H^{ev}(X_{g_0},\CC)
\end{equation*}
of the differential form
\begin{equation*}
 \ch_\Gamma p:=\tau \bl(e^{-\Omega/2\pi i}\br)\in
\bigoplus_{\left<g_0\right>\subset
\Gamma}\Lambda^{ev}(X_{g_0}).
\end{equation*}
\end{definition}

This is well defined. More precisely, the form $\ch_\Gamma
p$ is closed, and its cohomology class is independent of
the choice of a connection in the bundle $E$ and is
uniquely determined by the class of the projection $p$ in
the $K$-group $K_0(C^\infty(X,\End E)_\Gamma)$. The proof
is based on the identity
\begin{equation*}
 d\tau(A)=\tau\left([\nabla,A]\right),
\end{equation*}
where $A\in \Lambda^*(X,\End E)_\Gamma$ is an arbitrary
element such that $pA=A=Ap$.

\paragraph{Chern character of the symbol.}
Now let
\begin{equation*}
 D\colon C^\infty(M,E)\longrightarrow C^\infty(M,F)
\end{equation*}
be an elliptic operator with shifts acting in sections of
vector bundles on $M$. To define the Chern character of the
symbol $\sigma(D)$, we introduce a projection and hence an
element in $K$-theory associated with the symbol. To this
end, we make use of the bundle
$$
2B^*M=S(T^*M\oplus 1)
$$
of unit spheres in the vector bundle $T^*M\oplus 1$ over
$M$.

Consider the projection $p$ over the algebra $C(2B^*M,
\End(E\oplus F))_\Gamma$ defined by the formula
\begin{equation}\label{projector1}
p(\xi\cos\psi,\sin\psi)=\frac12
\begin{pmatrix}
  (1+\sin\psi)\id_E  &  \sigma^{-1}(D)(\xi)\cos\psi \\
  \sigma(D)(\xi)\cos\psi  & (1-\sin\psi)\id_F
\end{pmatrix},
\end{equation}
where $\xi$ lies on the unit sphere in $T^*M$, so that
$(\xi\cos\psi,\sin\psi)$ just lies on the unit sphere in
$T^*M\oplus 1$.

\begin{remark}
 Note that in general the projection $p$ is only continuous
but not infinitely differentiable at the points where
$\cos\psi=0$.
\end{remark}

We set
$$
[\sigma(D)]\ovs{\op{def}}=[p]\in
K_0\Bigl(C(2B^*M,\End(E\oplus F))_\Gamma\Bigr).
$$
Note that
\begin{equation*}
    K_0\bl(C(2B^*M,\End(E\oplus
F))_\Gamma\br)= K_0\bl(C^\infty(2B^*M,\End(E\oplus
F))_\Gamma\br),
\end{equation*}
since, as was already mentioned above,
$C^\infty(2B^*M,\End(E\oplus F))_\Gamma$ is a dense local
subalgebra of $C(2B^*M,\End(E\oplus F))_\Gamma$. Hence we
obtain the cohomology class
$$
\ch_\Gamma[\sigma(D)]\in \bigoplus_{g_0}
H^{ev}(2B^*M_{g_0},\mathbb{C}),
$$
which will be called the \textit{Chern character} of the
symbol $\si(D)$.

\subsection{Index theorem}

Now we are in position to state our main result.

\begin{theorem}
Let $D$ be an elliptic operator with shifts on the manifold
$M$. Then the index of $D$ is given by the formula
\begin{equation}
\label{cohom1} \ind D=\left\langle
\ch_\Gamma[\sigma(D)]\Td(TM\otimes\mathbb{C};\Gamma),
[2B^*M;\Gamma]\right\rangle,
\end{equation}
where
\begin{equation*}
 [2B^*M;\Gamma]=\prod_{g_0}[2B^*M_{g_0}]\in
\prod_{g_0} H_{ev}(2B^*M_{g_0})
\end{equation*}
is the fundamental class and angle brackets denote the
natural paring between cohomology and homology.
\end{theorem}

The proof involves extensive computations and goes by
reduction to the local index formula obtained for elliptic
operators with shifts in the first part of this paper.


\begin{thebibliography}{1}

\bibitem{0706.3511}
V.~E. Nazaikinskii, A.~Yu Savin, and B.~Yu Sternin.
\newblock On elliptic differential operators with shifts, 2007.
\newblock \url{http://arxiv.org/abs/0706.3511}.

\bibitem{AtSi3}
M.~F. Atiyah and I.~M. Singer.
\newblock The index of elliptic operators {III}.
\newblock {\em Ann. Math.}, {\bf 87}, 1968,  546--604.

\bibitem{BaCo2}
Paul Baum and Alain Connes.
\newblock Chern character for discrete groups.
\newblock In {\em A f\^ete of topology}, 1988, pages 163--232. Academic Press,
  Boston, MA.

\end{thebibliography}
\end{document}